\documentclass[preprint,1p,times,sort&compress]{elsarticle}
\usepackage{amsmath,amsfonts,mathrsfs}
\usepackage{amssymb}
\usepackage{color}
\usepackage{mathrsfs}
\usepackage{mathptmx}
\usepackage{verbatim}
\usepackage{epsfig}
\usepackage{CJK}
\usepackage{bm}
\usepackage{amsfonts}
\usepackage{amsthm}

\usepackage{amsthm}
\newtheorem{theorem}{Theorem}[section]

\newdefinition{rem}{Remark}

\numberwithin{equation}{section}

\journal{}

\begin{document}
\begin{frontmatter}
\title{Boundedness of solutions for Duffing equation with low regularity in time}

\author[SDUW]{Xiaoping Yuan\corref{cor}
}
\cortext[cor]{Supported by NSFC11421061.} \ead{xpyuan@fudan.edu.cn}
\address[SDUW]{School of Mathematical Sciences, Fudan University, Shanghai 200433, China}

\begin{abstract}
It is shown that all solutions are bounded for Duffing equation $\ddot{x}+ x^{2n+1}+\sum_{j=0}^{2n}P_{j}(t)x^{j}=0,$  provided that for each $n+1\le j\le 2n$, $P_j(t)\in C^{\gamma}(\mathbb T)$ with $\gamma>1-\frac1n$ and for each $0\le j\le n$, $P_j\in L(\mathbb T^1)$.
\end{abstract}


\end{frontmatter}

\section{Introduction}

In 1962, Moser \cite{Moser1962} proposed to study the boundedness of all solutions (Lagrange stability) for
Duffing equation
$$\ddot{x}+\beta x^{3}+\alpha x=P(t),\;\;P\in C(\mathbb{T}^{1}),\;\;\mathbb{T}^{1}:=\mathbb{R}/\mathbb{Z},$$
where $\beta>0,$ $\alpha\in \mathbb{R}$ are constants.

In 1976, Morris \cite{Morris1976} proved the boundedness of all solutions for
$$\ddot{x}+2 x^{3}=P(t).$$
Subsequently, Morris' boundedness results was, by Dieckerhoff-Zehnder \cite{Dieckerhoff-Zehnder1987} in 1987, extended to a wider class of systems
\begin{equation}\label{eq1}
\ddot{x}+ x^{2n+1}+\sum_{j=0}^{2n}P_{j}(t)x^{j}=0,\;n\geq 1,\end{equation}
where
$$P_{j}\in C^{\nu},\;\;\nu\geq 1+\frac{4}{n}+[\log_2^{n}]\rightarrow \infty, \;\;
\nLeftarrow\mbox{as}\;\; n \to \infty.$$
Then they remarked that

\ \

{\it
"It is not clear whether the boundedness phenomenon is related to the smoothness in the $t$-variable or whether this requirement is a shortcoming of our proof."}

\ \

In 1989 and 1992, Liu \cite{Liu1989, Liu1992} proved the boundedness for
$$\ddot{x}+ x^{2n+1}+a(t)x+P(t)=0,\;\;a(t)\in C^{0}(\mathbb{T}^{1}),\;\;P(t)\in C^{0}(\mathbb{T}^{1}).$$

In 1991, Laederich-Levi \cite{Laederich-Levi1991} relaxed the smoothness requirement of $P_{j}(t)\;(j=0, 1, \cdots, 2n)$ for \eqref{eq1} to
$$P_{j}\in C^{5+\varepsilon}(\mathbb T^1),\;\;\varepsilon>0.$$

In his  PhD thesis (1995), the present author further relaxed the requirement to $C^{2}$. See \cite{Yuan1995},\cite{Yuan1998} and \cite{Yuan2000}.

In the present paper, we will relax the smoothness requirement to $C^{1}.$ More exactly, we have the following theorem

\begin{theorem}\label{thm} For Arbitrary given constant $\gamma>1-\frac1n$, assume $P_{j}(t)\in C^{\gamma}(\mathbb{T}^{1})$ for $n+1\le j\le 2n$ and $P_j(t)\in L(\mathbb T^1)$ for $0\le j\le n$. Then every solution $x(t)$ of the equation \eqref{eq1},
\[\ddot{x}+ x^{2n+1}+\sum_{j=0}^{2n}P_{j}(t)x^{j}=0,\;n\geq 1,\]
 is bounded, i.e. it exists for all $t\in \mathbb{R}$ and
$$\sup_{t\in \mathbb{R}}(|x(t)|+|\dot{x}(t)|)<C<\infty,$$
where $C=C(x(0), \dot{x}(0))$ depends the initial data $(x(0), \dot{x}(0)).$
\end{theorem}
 \begin{rem} In \cite{Wang2000}, it is proved that there is  a continuous periodic function $p(t)$  such that the Duffing equation $\frac{d^2x}{dt^2}+x^{2n+1}+p(t) x^l=0$ with $p(t)\in C^0(\mathbb T^1)$, $n\ge2$, $2n+1>l\ge n+2$ possesses an unbounded solution, which shows that the smoothness of the coefficients $P_j(t)$'s  does influence the boundedness of solutions. Therefor, the result of theorem \ref{thm} is sharp without considering the derivative of non-integral order.

 \end{rem}

 \section{Action-Angle Variable}
 Replacing $x$ by $Ax$ in \eqref{eq1}, we get
 \begin{equation}\label{eq2}
 A\ddot{x}+A^{2n+1}x^{2n+1}+\sum_{j=0}^{2n}P_{j}(t)x^{j}A^{j}=0,
 \end{equation}
 where $A$ is a constant large enough.
 That is,
 \begin{equation}\label{eq3}
 \ddot{x}+A^{2n}x^{2n+1}+\sum_{j=0}^{2n}P_{j}(t)x^{j}A^{j-1}=0.
 \end{equation}
 Let
 $$y=A^{-n}\dot{x},\;\;\mbox{or}\;\;\dot{x}=A^{n}y.$$
 Then
 \begin{eqnarray*}
 \dot{y}&=&A^{-n}\ddot{x}\\
 &=& A^{-n}(-A^{2n}x^{2n+1}-\sum_{j=0}^{2n}P_{j}(t)x^{j}A^{j-1}\\
 &=&-A^{n}x^{2n+1}-\sum_{j=0}^{2n}P_{j}(t)x^{j}A^{j-n-1}.
 \end{eqnarray*}
 Thus,
 \begin{equation}\label{eq4}
 \dot{x}=\frac{\partial H}{\partial y},\;\;\dot{y}=-\frac{\partial H}{\partial x},
 \end{equation}
 where
 \begin{equation}\label{eq5}
 H=A^{n}\left(\frac{1}{2}y^{2}+\frac{1}{2(n+1)}x^{2(n+1)}\right)+\sum_{j=0}^{2n}\frac{P_{j}(t)}{j+1}x^{j+1}A^{j-n-1}.
 \end{equation}
 Let
 $$\mathbb{T}^{1}_{s}=\{t\in \mathbb{C}/\mathbb{Z}: |\mathrm{Im} \;t|<s\} \;\;\mbox{for\;any}\;\;s>0.$$
 Consider an auxiliary Hamiltonian system
 \begin{equation}\label{eq6}
 \dot{x}=\frac{\partial H_{0}}{\partial y},\;\;\dot{y}=-\frac{\partial H_{0}}{\partial x},\;\;
 H_{0}=\frac{1}{2}y^{2}+\frac{1}{2(n+1)}x^{2(n+1)}.
 \end{equation}
Let $(x_{0}(t), y_{0}(t))$ be the solution to \eqref{eq6} with initial $(x_{0}(0), y_{0}(0))=(1, 0).$
Then this solution is clearly periodic. Let $T_{0}$ be its minimal positive period. By Energy conservation,
we has
\begin{equation}\label{eq7}
(n+1)y_{0}^{2}(t)+x_{0}^{2n+2}(t)\equiv 1,\;\;t\in \mathbb{R},
\end{equation}
 by which, we construct the following symplectic transformation
$$\Psi_{0}:\;\;
\left\{\begin{array}{ll}
x=c^{\alpha}I^{\alpha}x_{0}(\theta T_{0}),\\
\\
y=c^{\beta}I^{\beta}y_{0}(\theta T_{0}),
\end{array}
\right.$$
where $\alpha=\frac{1}{n+2},\;\beta=1-\alpha=\frac{n+1}{n+2},\; c=\frac{1}{\alpha T_{0}}$
and where $(I, \theta)\in \mathbb{R}^{+}\times \mathbb{T}^{1}$ is action-angle variables.
By calculation, $\text{det}\, \frac{\partial\,(x,y)}{\partial(\theta,I)}=1$. Thus the transformation is indeed symplectic.
  %
 Clearly $\Psi_{0}(I, \theta)$ is analytic in $(I, \theta)\in \mathbb{R}^{+}\times \mathbb{T}^{1}_{s_{0}}$ with some constant $s_{0}>0.$

Under $\Psi_{0},$ \eqref{eq4} is changed
\begin{equation}\label{eq8}
\dot{\theta}=\frac{\partial H}{\partial I},\;\;\dot{I}=-\frac{\partial H}{\partial \theta},
\end{equation}
where $H=H_{0}(I)+R(I, \theta, t)$ with
\begin{equation}\label{eq9}
H_{0}(I)=d\cdot A^{n}\cdot I^{2\beta}=d \cdot A^{n}\cdot I^{\frac{2(n+1)}{n+2}},\;\;d=\frac{c^{2\beta}}{2(n+1)},
\end{equation}
and
\begin{equation}\label{eq10}
R(I, \theta, t)=\sum_{j=0}^{2n}\frac{P_{j}(t)}{j+1}(c^{\frac{1}{n+1}}\,x_{0}(\theta T_{0}))^{j+1}A^{j-n-1}I^{\frac{j+1}{n+2}}.
\end{equation}
Clearly, $R(I, \theta, t)=O(A^{n-1})$ for $A \to \infty$ and fixed $I$ belongs to some compact intervals.
\section{Approximation Lemma}
First, we cite an approximation lemma.
See \cite{Salamon1989} and \cite{Salamon2004}, for the detail.
We start by recalling some definitions and setting some new notations. Assume $X$ is a Banach space with the norm
$||\cdot||_{X}$. First recall that $C^{\mu}(\mathbb{R}^{n}; X)$ for $0< \mu <1$ denotes the space of bounded
H\"{o}lder continuous functions $f: \mathbb{R}^{n}\mapsto X$ with the form
$$\|f\|_{C^{\mu}, X}=\sup_{0<|x-y|<1}\frac{\|f(x)-f(y)\|_{X}}{|x-y|^{\mu}}+\sup_{x\in \mathbb{T}^{n}}\|f(x)\|_{X}.$$
If $\mu=0$ then $\|f\|_{C^{\mu},X}$ denotes the sup-norm. For $\ell=k+\mu$ with $k\in \mathbb{N}$ and $0\leq \mu <1,$
we denote by $C^{\ell}(\mathbb{R}^{n};X)$ the space of functions $f:\mathbb{R}^{n}\mapsto X$ with H\"{o}lder continuous partial derivatives, i.e., $\partial ^{\alpha}f\in C^{\mu}(\mathbb{R}^{n}; X_{\alpha})$ for all muti-indices $\alpha=(\alpha_{1}, \cdots, \alpha_{n})\in \mathbb{N}^{n}$ with the assumption that
$|\alpha|:=|\alpha_{1}|+\cdots+|\alpha_{n}|\leq k$ and $X_{\alpha}$ is the Banach space of bounded operators
$T:\prod^{|\alpha|}(\mathbb{R}^{n})\mapsto X$ with the norm
$$\|T\|_{X_{\alpha}}=\sup \{||T(u_{1}, u_{2}, \cdots, u_{|\alpha|})||_{X}:\|u_{i}\|=1, \;1\leq i \leq |\alpha|\}.$$
We define the norm
$$||f||_{C^{\ell}}=\sup_{|\alpha|\leq \ell}||\partial ^{\alpha}f||_{C^{\mu}, X_{\alpha}}$$

\begin{theorem}(Jackson)\label{thm3.1}
Let $f\in C^{\ell}(\mathbb{R}^{n}; X)$ for some $\ell>0$ with finite $C^{\ell}$ norm over $\mathbb{R}^{n}.$
Let $\phi$ be a radical-symmetric, $C^{\infty}$ function, having as support the closure of the unit ball centered at the origin, where $\phi$ is completely flat and takes value 1, let $K=\widehat{\phi}$ be its Fourier transform. For all $\sigma >0$ define
$$ f_{\sigma}(x):=K_{\sigma}\ast f=\frac{1}{\sigma^{n}}\int_{\mathbb{T}^{n}}K(\frac{x-y}{\sigma})f(y)dy.$$
Then there exists a constant $C\geq 1$ depending only on $\ell$ and $n$ such that the following holds: For any $\sigma >0,$ the function $f_{\sigma}(x)$ is a real-analytic function from $\mathbb{C}^{n}$ to $X$ such that if $\Delta_{\sigma}^{n}$ denotes the $n$-dimensional complex strip of width $\sigma,$
$$\Delta_{\sigma}^{n}:=\{x\in \mathbb{C}^{n}\big ||\mathrm{Im} x_{j}|\leq \sigma,\;1\leq j\leq n\},$$
then for $\forall \alpha\in\mathbb{N}^{n}$ with $|\alpha|\leq \ell$ one has

\begin{equation}\label{cite3.1}
\sup_{x\in \Delta_{\sigma}^{n}}||\partial ^{\alpha}f_{\sigma}(x)
-\sum_{|\beta|\leq \ell-|\alpha|}\frac{\partial^{\beta+\alpha}f(\mathrm{Re}x)}{\beta !}(\sqrt{-1}\mathrm{Im}x)^{\beta}||_{X_{\alpha}}\leq C ||f||_{C^{\ell}}\sigma^{\ell-|\alpha|},
\end{equation}

and for all $0\leq s\leq \sigma,$
\begin{equation}\label{cite3.2}
\sup_{x\in \Delta_{s}^{n}}\|\partial^{\alpha} f_{\sigma}(x)-\partial^{\alpha}f_{s}(x)\|_{X_{\alpha}}
\leq C ||f||_{C^{\ell}}\sigma^{\ell-|\alpha|}.
\end{equation}

The function $f_{\sigma}$ preserves periodicity (i.e., if $f$ is T-periodic in any of its variable $x_{j}$, so is $f_{\sigma}$).
\end{theorem}
%
  %
  By this theorem, for each $P_{j}(t)\in C^{\gamma}(\mathbb{T}^{1}),\;j=n+1, 1, \cdots, 2n,$ and any $\varepsilon>0,$
there is a real analytic function \footnote{A complex value function $f(t)$ of complex variable $t$ in some domain in $\mathbb C$ is called real analytic if it is analytic in the domain and is real for real argument $t$} $\,$ $P_{j,\varepsilon}(t)$ from $\mathbb{T}^{1}_{\varepsilon}$ to $\mathbb{C}$ such
that
\begin{equation}\label{eq11}
\sup_{t\in\mathbb{T}^{1}}|P_{j,\varepsilon}(t)-P_{j}(t)|\leq C \,\varepsilon^{\gamma}\, ||P_{j}||_{C^{\gamma}},
\end{equation}
and
\begin{equation}\label{eq12}
\sup_{t\in\mathbb{T}_{\varepsilon}^{1}}|P_{j,\varepsilon}(t)|\leq C ||P_{j}||_{C^{\gamma}}.
\end{equation}

Write
\begin{equation}\label{eq13}
R(I, \theta, t)=R_{\varepsilon}(I, \theta, t)+R^{\varepsilon}(I, \theta, t),
\end{equation}
where
\begin{equation}\label{eq14}
R_{\varepsilon}(I, \theta, t)=\sum_{j=n+1}^{2n}\frac{1}{j+1}A^{j-n-1}I^{\frac{j+1}{n+2}}
c^{\frac{j}{n+2}}x_{0}^{j+1}(\theta T_{0})P_{j, \varepsilon}(t),
\end{equation}

\begin{equation}\label{eq15}\begin{array}{ll}
R^{\varepsilon}(I, \theta, t)=&\sum_{j=0}^{n}\frac{1}{j+1}A^{j-n-1}I^{\frac{j+1}{n+2}}
c^{\frac{j}{n+2}}x_{0}^{j+1}(\theta T_{0})P_{j}(t)\\&+\sum_{j=n+1}^{2n}\frac{1}{j+1}A^{j-n-1}I^{\frac{j+1}{n+2}}
c^{\frac{j}{n+2}}x_{0}^{j+1}(\theta T_{0})(P_{j}(t)-P_{j, \varepsilon}(t)).\end{array}
\end{equation}
Now let us restrict $I$ belongs to some compact intervals, $[1, 4],$ say. Let
\[A^{-1}<\varepsilon_0.\]

For a sufficiently small $\varepsilon_0>0$, letting
\begin{equation}\label{eq16}
\varepsilon=(\varepsilon_{0}/A^{n-1})^{1/\gamma},
\end{equation}
 by Theorem \ref{thm3.1}, we have the following facts:
\begin{description}
       \item[(i)]$R^{\varepsilon}(I, \theta, t)$ is real analytic in $(I, \theta)\in [1, 4]\times \mathbb{T}^{1}_{s_0}$ for
fixed $t\in \mathbb{T}^{1}$ and $R^{\varepsilon}(I,\theta,\cdot)\in L^1(\mathbb{T}^{1})$ for fixed $(I, \theta)\in [1, 4]\times \mathbb{T}^{1}_{s_{0}},$ and
\begin{equation}\label{eq17}
\sup_{(I, \theta, t)\in [1, 4]\times \mathbb{T}^{1}_{s_{0}}\times \mathbb{T}^{1}}|R^{\varepsilon}(I, \theta, t)|\leq C\varepsilon_{0},
\end{equation}
where $C$ is a constant \footnote{Denote by $C$ a universal constant which may be different in different place.} depending on only $||P_{j}||_{C^{\gamma}}$.
       \item[(ii)] $R_{\varepsilon}(I, \theta, t)$ is real analytic in $(I, \theta, t)\in[1, 4]\times \mathbb{T}^{1}_{s_{0}}\times \mathbb{T}^{1}_{\varepsilon}$ and
\begin{equation}\label{eq18}
\sup_{(I, \theta, t)\in [1, 4]\times \mathbb{T}^{1}_{s_{0}}\times \mathbb{T}^{1}_{\varepsilon}}
|R_{\varepsilon}(I, \theta, t)|\leq CA^{n-1},
\end{equation}
where $C$ is a constant depends on only $||P_{j}||_{C^{\gamma}}.$
Therefore, we have
\begin{equation}\label{eq19}
H(I, \theta, t)=H_{0}(I)+R_{\varepsilon}(I, \theta, t)+R^{\varepsilon}(I, \theta, t).
\end{equation}
\end{description}
\section{Symplectic transformations}
We will look for a series of symplectic transformations $\Psi_{1}, \cdots, \Psi_{N}$
such that
\begin{equation*}
H^{(N)}=H\circ \Psi_{1}\circ \cdots\circ \Psi_{N}=H_{0}^{N}+O(\varepsilon_{0}),
\end{equation*}
where $H_{0}^{N}(\mu)\approx A^{n}\mu^{\frac{2(n+1)}{n+2}}$
and that Moser's twist works for $H^{(N)}.$

To this end, let $\Psi_{1}: (\mu, \phi)\mapsto (I, \theta)$ is implicitly defined by
$$\Psi_{1}: \left\{
              \begin{array}{ll}
                I=\mu+\frac{\partial S_{1}}{\partial \theta} \\
               \phi=\theta+\frac{\partial S_{1}}{\partial \mu}
              \end{array}
            \right.
$$
with $S_{1}=S_{1}(\mu, \theta, t)$ to be specified latter. If $\Psi_{1}$ is well-defined, then it is symplectic,
since
$$dI \wedge d\theta=(1+\frac{\partial^{2}S_{1}}{\partial \mu \partial \theta})d\mu \wedge d\theta=d\mu \wedge d \phi.$$
The transformed Hamiltonian function $H^{(1)}(\mu, \phi, t)=H\circ \Psi_{1}(\mu, \phi, t).$ We express temporarily in the variable $(\mu, \theta)$ instead of $(\mu, \phi)$:
\begin{equation}\label{eq40}
H^{(1)}(\mu, \theta, t)=H(\mu+\frac{\partial S_{1}}{\partial \theta}, \theta, t)+\frac{\partial S_{1}}{\partial t}.
\end{equation}
By Taylor's formula and \eqref{eq19}
\begin{eqnarray}\label{eq41}
{H}^{(1)}(\mu, \theta, t)&=&H_{0}(\mu+\frac{\partial S_{1}}{\partial \theta}, \theta, t)+R_{\varepsilon}(\mu+\frac{\partial S_{1}}{\partial \theta}, \theta, t)+R^{\varepsilon}\circ \Psi_{1}(\mu, \phi, t)+\frac{\partial S_{1}}{\partial t}\nonumber\\
&=& H_{0}(\mu)+\partial_{\mu}H_{0}(\mu)\frac{\partial S_{1}}{\partial \theta}+R_{\varepsilon}(\mu, \theta, t)+
R^{1}_{\varepsilon}(\mu, \theta, t)+R^{\varepsilon}\circ \Psi_{1}(\mu, \phi, t),
\end{eqnarray}
where
\begin{eqnarray}\label{eq42}
R^{1}_{\varepsilon}(\mu, \theta, t)&=&\int_{0}^{1}(1-\tau)\partial ^{2}_{\mu}H_{0}(\mu+\frac{\partial S_{1}}{\partial \theta}\tau, \theta, t)(\frac{\partial S_{1}}{\partial \theta})^{2}d\tau\nonumber\\
&&+\int_{0}^{1}\partial_{\mu}R_{\varepsilon}(\mu+\frac{\partial S_{1}}{\partial \theta}\tau, \theta, t)\frac{\partial S_{1}}{\partial \theta}d\tau+\frac{\partial S_{1}}{\partial t}.
\end{eqnarray}
Let
\begin{equation}\label{eq43}
\partial_{\mu}H_{0}\cdot \frac{\partial S_{1}}{\partial \theta}+R_{\varepsilon}(\mu, \theta, t)=[R_{\varepsilon}](\mu, t),\;\;
[R_{\varepsilon}](\mu, t)=\int_{0}^{1}R_{\varepsilon}(\mu, \theta, t)d\theta.
\end{equation}
Then
\begin{eqnarray}\label{eq44}
H^{(1)}(\mu, \theta, t)&=&H_{0}(\mu)+[R_{\varepsilon}](\mu, t)+R^{1}_{\varepsilon}(\mu, \theta, t)+R^{\varepsilon}\circ \Psi_{1}(\mu, \phi, t)\nonumber\\
&=&H^{1}_{0}(\mu, t)+R^{1}_{\varepsilon}(\mu, \theta, t)+R^{\varepsilon}\circ \Psi_{1}(\mu, \phi, t),
\end{eqnarray}
where
\begin{equation}\label{eq45}
H^{1}_{0}(\mu, t)=H_{0}(\mu)+[R_{\varepsilon}](\mu, t).
\end{equation}
We are now in position to solve \eqref{eq43}.
\begin{equation}\label{eq46}
S_{1}(\mu, \theta, t)=\int_{0}^{\theta}\frac{[R_{\varepsilon}](\mu, t)-R_{\varepsilon}(\mu, \theta, t)}{\partial_{\mu}H_{0}(\mu)}d\theta.
\end{equation}
By \eqref{eq16} and \eqref{eq18}, $S_{1}$ is well-defined in $(\mu, \theta, t)\in [1, 4]\times \mathbb{T}_{s_{0}}\times \mathbb{T}_{\varepsilon},$ and analytic in the domain, and
\begin{equation}\label{eq47}
\sup_{(\mu, \theta, t)\in [1, 4]\times \mathbb{T}_{s_{0}}\times \mathbb{T}_{\varepsilon}}|S_{1}(\mu, \theta, t)|
\leq CA^{-1}.
\end{equation}
Thus, by the implicit function theorem, $\Psi_{1}(\mu, \phi, t): [1+O(A^{-1}), 4-O(A^{-1})]\times \mathbb{T}^1_{s_{0}/2}\times \mathbb{T}^1_{\varepsilon}\to [1, 4]\times\mathbb{T}^1_{s_{0}}\times \mathbb{T}^1_{\varepsilon}.$
\begin{itemize}
  \item Estimate of $H^{1}_{0}(\mu, t).$

By \eqref{eq18}, we have that $H^{1}_{0}(\mu, t)$ is analytic in $[1, 4]\times \mathbb{T}_{\varepsilon},$ and
\begin{equation}\label{eq48}
CA^{n}\geq |\partial_{\mu}^2\,H_{0}^{1}(\mu, t)|\geq \frac{A^{n}}{C}, \;\;t\in \mathbb{T}_{\varepsilon/2},
\end{equation}
and by Cauchy's estimate
\begin{eqnarray}\label{eq49}
&&\sup_{(\mu, t)\in [1, 4]\times \mathbb{T}^1_{\varepsilon/2}}|\partial_{t}H_{0}^{1}(\mu, t)|\nonumber\\
&\leq &
\sup_{(\mu, t)\in [1, 4]\times \mathbb{T}^1_{\varepsilon/2}}|\partial_{t}[R_{\varepsilon}](\mu, t)|\nonumber\\
&\leq &\sup_{(\mu, \theta, t)\in [1, 4]\times \mathbb{T}^1_{s_{0}}\times \mathbb{T}^1_{\varepsilon/2}}
|\partial_{t}R_{\varepsilon}(\mu, \theta, t)|\nonumber\\
&\leq & \frac{2}{\varepsilon}\sup_{(\mu, \theta, t)\in [1, 4]\times \mathbb{T}^1_{s_{0}}\times \mathbb{T}^1_{\varepsilon}}
|R_{\varepsilon}(\mu, \theta, t)|\nonumber\\
&\leq & \frac{2}{\varepsilon} CA^{n-1}\lesssim C\varepsilon_{0}^{-\frac{1}{\gamma}}A^{\frac{n-1}{\gamma}}A^{n-1}\nonumber\\
&\leq &  C\varepsilon_{0}^{-\frac{1}{\gamma}}A^{(n-1)(1+\frac{1}{\gamma})}.
\end{eqnarray}
  \item Estimate of $R^{1}_{\varepsilon}(\mu, \theta, t).$

By \eqref{eq47} and the Cauchy estimate,
\begin{equation}\label{yuan1}\begin{array}{ll}\sup_{(\mu, \theta, t)\in [1, 4]\times \mathbb{T}_{s_0}\times \mathbb{T}_{\varepsilon/2}}|\partial_{t}S_{1}(\mu, \theta, t)|&\leq \frac{2CA^{-1}}{\varepsilon}\leq CA^{-1}(\frac{\varepsilon_{0}}{A^{n-1}})^{-\frac{1}{\gamma}}
\leq C\varepsilon_{0}^{-\frac{1}{\gamma}}A^{-1+\frac{n-1}{\gamma}}\\&
=C\varepsilon_{0}^{-\frac{1}{\gamma}}A^{n-1-\varpi}.\end{array}\end{equation}
where \begin{equation}\label{yuan2}\varpi:=n-\frac{n-1}{\gamma}=\frac{n}{\gamma}\left(\gamma-(1-\frac{1}{n})
\right). \end{equation}
By assuming $ 1\ge\gamma>1-\frac1n$,
\[0<\varpi\le 1.\]

By \eqref{eq18} and noting $H_{0}(\mu)=d A^{n-1}\mu^{\frac{2n+2}{n+2}},$ we have
\begin{eqnarray}\label{eq50}
\sup_{(\mu, \theta, t)\in \mathcal{D}_{1}}|R^{1}_{\varepsilon}(\mu, \theta, t)|
&\leq &CA^{n}A^{-2}+CA^{n-1}A^{-1}+C\varepsilon_{0}^{-\frac{1}{\gamma}}A^{-1+\frac{n-1}{\gamma}}\nonumber\\
&\leq & C\,\varepsilon_{0}^{-\frac{1}{\gamma}}\,A^{n-1-\varpi},
\end{eqnarray}
where
$$\mathcal{D}_{1}=[1+O(A^{-1}), 4-O(A^{-1})]\times \mathbb{T}^1_{s_{0}/2}\times \mathbb{T}^1_{\varepsilon/2},\;\;
.$$
By \eqref{eq47} and the implicit function theorem, there exist $U_{1}(\mu, \phi, t),\;V_{1}(\mu, \phi, t)$ analytic in $\mathcal{D}_{1}$ such that
\begin{equation}\label{eq51}
\sup_{\mathcal{D}_{1}}|U_{1}|\leq CA^{-1},\;\;\sup_{\mathcal{D}_{1}}|V_{1}|\leq CA^{-1},
\end{equation}
and
\begin{equation}\label{eq52}
\Psi_{1}: \left\{
            \begin{array}{ll}
              I=\mu+U_{1}(\mu, \phi, t) \\
              \theta=\phi+V_{1}(\mu, \phi, t)
            \end{array}
          \right.
\end{equation}
and
\begin{equation}\label{eq53}
H^{1}(\mu, \phi, t)=H^{1}_{0}(\mu, t)+\widetilde{R}^{1}_{\varepsilon}(\mu, \phi, t)
+R^{\varepsilon}\circ \Psi(\mu, \phi, t),
\end{equation}
where
\begin{equation}\label{eq54}
\widetilde{R}^{1}_{\varepsilon}(\mu, \phi, t)={R}^{1}_{\varepsilon}(\mu, \phi+V_{1}(\mu, \phi, t), t)
\end{equation}
and
\begin{equation}\label{eq55}
\sup_{\mathcal{D}_{1}}|\widetilde{R}^{1}_{\varepsilon}(\mu, \phi, t)|\leq C\,\varepsilon_{0}^{-\frac{1}{\gamma}}\,A^{n-1-\varpi}.
\end{equation}
Similarly, let
\begin{equation}\label{eq56}
\Psi_{2}: \left\{
            \begin{array}{ll}
              \mu=\lambda+\frac{\partial S_{2}}{\partial \phi}\\
              \widetilde{\phi}=\phi+\frac{\partial S_{2}}{\partial \lambda}
            \end{array}
          \right.
\end{equation}
with $S_{2}=S_{2}(\lambda, \phi, t)$ is defined by
\begin{equation}\label{eq57}
S_{2}(\lambda, \phi, t)=\int_{0}^{\phi}\frac{[\tilde R_{\varepsilon}^{1}](\lambda, t)-\tilde R^{1}_{\varepsilon}(\lambda, \phi, t)}{\partial_{\mu}H_{0}^{1}(\mu, t)}dt,\quad [\tilde R_{\varepsilon}^{1}](\lambda, t)=\int_0^1 \, \tilde R^{1}_{\varepsilon}(\lambda, \phi, t)\, d\phi.
\end{equation}
By \eqref{eq48} and \eqref{eq50},
\begin{equation}\label{eq58}
\sup_{\mathcal{D}_{1}}|S_{2}(\lambda, \phi, t)|\leq C\,A^{-n}\,(\frac{1}{\varepsilon_{0}})^{\frac{1}{\gamma}}\,A^{n-1-\varpi}\leq C\,(\frac{1}{\varepsilon_{0}})^{\frac{1}{\gamma}}\,A^{-1-\varpi}.
\end{equation}
It follows from the implicit function theorem that $\Psi_{2}: (\lambda, \widetilde{\phi})\in \mathcal{D}_{2}\to \mathcal{D}_{1}$ is well-defined, where $\mathcal{D}_{2}=[1+O(A^{-1}), 4-O(A^{-1})]\times \mathbb{T}_{s_{0}/4}\times \mathbb{T}_{\varepsilon/4}.$ By Cauchy estimate,
\begin{equation}\label{eq59}
\sup_{\mathcal{D}_{2}}|\partial_{t}S_{2}(\lambda, \phi, t)|\leq \frac{C}{\varepsilon}\,A^{-1-\varpi}\,(\frac{1}{\varepsilon_{0}})^{\frac{1}{\gamma}}
\leq C\,A^{-\varpi-1}(\frac{1}{\varepsilon_{0}})^{\frac{2}{\gamma}}A^{\frac{n-1}{\gamma}}=
C\,(\frac{1}{\varepsilon_{0}})^{\frac{2}{\gamma}}A^{n-2\varpi-1}.
\end{equation}
Let
\begin{eqnarray}\label{eq60}
&&H^{(2)}(\lambda, \phi, t):=H^{1}(\lambda+\frac{\partial S_{2}}{\partial \phi}, \phi, t)+\frac{\partial S_{2}}{\partial t}\nonumber\\
&=& H^{1}_{0}(\lambda, t)+\partial_{\lambda}H^{1}_{0}(\lambda, t)\frac{\partial S_{2}}{\partial \phi} +\tilde
R^{1}_{\varepsilon}(\lambda, \phi, t)+
R^{2}_{\varepsilon}(\lambda, \phi, t)\nonumber\\
&&+R^{\varepsilon}\circ \Psi_{1}\circ \Psi_{2}(\lambda, \widetilde{\phi}, t),
\end{eqnarray}
where
\begin{eqnarray}\label{eq61}
R^{2}_{\varepsilon}(\lambda, \phi, t)&=&\int_{0}^{1}(1-\tau)\partial_{\lambda}^{2}H_{0}^{1}(\lambda+\frac{\partial S_{2}}{\partial \phi}\tau, \phi, t)(\frac{\partial S_{2}}{\partial \phi})^{2}d\tau\nonumber\\
&&+\int_{0}^{1}\partial_{\lambda}\tilde R^{1}_{\varepsilon}(\lambda+\frac{\partial S_{2}}{\partial \phi}, \phi, t)\frac{\partial S_{2}}{\partial \phi}d\tau+\frac{\partial S_{2}}{\partial t}.
\end{eqnarray}
By \eqref{eq57},
$$\partial_{\lambda}H_{0}^{1}(\lambda, t)\frac{\partial S_{2}}{\partial \phi}
+\tilde R_{\varepsilon}^{1}(\lambda, \phi, t)=[\tilde R^{1}_{\varepsilon}](\lambda, t).$$
Let
\begin{eqnarray}\label{eq62}
H^{2}_{0}(\lambda, t)=H_{0}^{1}(\lambda, t)+[\tilde R^{1}_{\varepsilon}](\lambda, t).
\end{eqnarray}
It follows that
\begin{eqnarray}\label{eq63}
H^{2}(\lambda, \phi, t)=H^{2}_{0}(\lambda, t)+R_{\varepsilon}^{2}(\lambda, \phi, t)+R^{\varepsilon}\circ \Psi_{1} \circ\Psi_{2}(\lambda, \widetilde{\phi}, t).
\end{eqnarray}
  \item Estimate of $H^{2}_{0}(\lambda, t).$

By \eqref{eq48}, \eqref{eq49} and \eqref{eq50},
\begin{eqnarray}\label{eq64}
CA^{n}\geq |\partial_{\lambda}^2H^{2}_{0}(\lambda, t)|\geq \frac{A^{n}}{C},\;\;\lambda\in [1, 4],\;\;t\in \mathbb{T}_{\varepsilon/2}.
\end{eqnarray}

\begin{eqnarray}\label{eq65}
\sup_{(\lambda,t)\in [1, 4]\times \mathbb{T}^1_{\varepsilon/4}}|\partial _{t}H^{2}_{0}(\lambda, t)|\leq C\varepsilon_{0}^{-\frac{1}{\gamma}}A^{(n-1)(1+\frac{1}{\gamma})}.
\end{eqnarray}
\item Estimate of $R^{2}_{\varepsilon}(\lambda, \phi, t).$

By \eqref{eq50}, \eqref{eq58}, \eqref{eq59} and \eqref{eq61}, \eqref{eq64},

\begin{eqnarray}\label{eq66-1}\begin{array}{ll}&
\sup_{\mathcal{D}_{2}}|R^{2}_{\varepsilon}(\lambda, \phi, t)|\\ &\leq CA^{n}(\frac{1}{\varepsilon_{0}})^{2/\gamma}A^{-2(1+\varpi)}+C (\frac{1}{\varepsilon_{0}})^{2/\gamma}A^{n-1-\varpi}\,A^{-1-\varpi}+
C\,(\frac{1}{\varepsilon_{0}})^{2/\gamma}A^{n-1-2\varpi}
\\& \leq C\,(\frac{1}{\varepsilon_{0}})^{2/\gamma}A^{n-1-2\varpi}.\end{array}
\end{eqnarray}
Take $N\in \mathbb{N}$ with $n-\varpi\,{N}\leq -1.$
Repeating the above procedure $N$ times, we get a series of symplectic transformations $\Psi_{1}, \cdots, \Psi_{N}$
such that
\begin{eqnarray*}
H^{N}(\rho, \xi, t)&=&H\circ\Psi_{1}\circ\cdots\circ\Psi_{N}\\
&=&H^{N}_{0}(\rho, t)+R^{N}_{\varepsilon}(\rho, \xi, t)+R^{\varepsilon}\circ\Psi_{1}\circ\Psi_{N}(\rho, \xi, t),
\end{eqnarray*}
where $(\rho, \xi, t)\in [1+O(A^{-1}), 4-O(A^{-1})]\times \mathbb{T}_{s_{0}/2^{N}}^1\times \mathbb{T}_{\varepsilon/2^{N}}^1,$
and
\begin{equation}\label{*}
\Phi\triangleq \Psi_{1}\circ \cdots \circ\Psi_{N}: [1+O(A^{-1}), 4-O(A^{-1})]\times \mathbb{T}^1\times \mathbb{T}^1\to[1, 4]\times \mathbb{T}^1\times \mathbb{T}^1,
\end{equation}
and
\begin{equation}\label{**}
\Phi=id+O(A^{-1}).
\end{equation}
and $H_{0}^{N}(\rho, t)$ satisfies
\begin{equation}\label{eq66}
CA^{n}\geq |\partial_{\rho}^2H_{0}^{N}(\rho, t)|\geq \frac{A^{n}}{C},\;\;\rho\in [2, 3],\;\;t\in \mathbb{T},
\end{equation}

\begin{equation}\label{eq67}
\sup_{(\rho, t)\in[2,3]\times \mathbb{T}}|\partial_{t}H_{0}^{N}(\rho, t)|\leq C\varepsilon_{0}^{-\frac{1}{\gamma}}A^{(n-1)(1+\frac{1}{\gamma})}
\end{equation}
and
$R_{\varepsilon}^{N}(\rho, \xi, t)$ satisfies that for $0\leq p+q \leq 6.$
\begin{equation}\label{eq68}
\sup_{(\rho, \xi, t)\in [2, 3]\times \mathbb{T}\times\mathbb{T}}|\partial^{p}_{\rho}\partial^{q}_{\xi}R^{N}_{\varepsilon}(\rho, \xi, t)|\leq C\,A^{n-\varpi{N}}\left(\frac{1}{\varepsilon_0}\right)^{N/\gamma}\leq C\,A^{-1}\left(\frac{1}{\varepsilon_0}\right)^{N/\gamma}<C\varepsilon_{0}
\end{equation} where $C$ depends on $N$ and  we have assumed that $A$ is large enough such that $$A^{-1}\left(\frac{1}{\varepsilon_0}\right)^{N/\gamma}<\varepsilon_0.$$
Let
\begin{equation}\label{eq69}
\mathcal{R}(\rho, \xi, t)=R^{N}_{\varepsilon}(\rho, \xi, t)+R^{\varepsilon}\circ\Psi^{1}\circ\cdots\circ\Psi^{N}.
\end{equation}
Then by \eqref{eq17}, \eqref{eq68}, \eqref{*} and \eqref{**}
\begin{equation}\label{eq70}
\sup_{(\rho, \xi)\in [2,3]\times \mathbb{T}^1}\int_0^1|\partial_{\rho}^{p}\partial_{\xi}^{q}\mathcal{R}(\rho, \xi, t)|\, dt\leq C\varepsilon_{0},\;\;0\leq p+q \leq 6.
\end{equation}
Now,
\begin{equation}\label{eq71}
H^{N}(p, \xi, t)=H_{0}^{N}(\rho, t)+\mathcal{R}(\rho, \xi, t).
\end{equation}
\end{itemize}
\section{Proof of theorem}
For $H^{N}$ the Hamiltonian equation is
\begin{equation}\label{eq72}
\left\{
  \begin{array}{ll}
    \dot{\rho}=-\frac{\partial H^{N}}{\partial \xi}=-\frac{\partial \mathcal{R}(\rho, \xi, t)}{\partial \xi} =O(\varepsilon_{0}),\\
 \dot{\xi}=\frac{\partial H^{N}}{\partial \rho}=\frac{\partial H^{N}_{0}(\rho, t)}{\partial \rho}+\frac{\partial \mathcal{R}(\rho, \xi, t)}{\partial \xi}
=\frac{\partial H^{N}_{0}(\rho, t)}{\partial \rho}+O(\varepsilon_{0}).
  \end{array}
\right.
\end{equation}
Note
$$H_{0}^{N}=d \cdot A^{n}\cdot\rho^{\frac{2(n+1)}{n+2}}+O(A^{n-1}).$$
By suing Picard iteration and Gronwall's inequality and noting \eqref{eq70}, we get that the time-1 map of \eqref{eq72} is of the form

$$\mathcal{P}: \left\{
                \begin{array}{ll}
                  \rho_{1}=\rho(t)|_{t=1}=\rho_{0}+{F}(\rho_{0}, \xi_{0}), \\
                  \xi_{1}=\xi(t)|_{t=1}=\xi_{0}+\alpha(\rho_0)+G(\rho_{0}, \xi_{0}),
                \end{array}
              \right.
(\rho_{0}, \xi_{0})\in [2,3]\times \mathbb{T}^1
$$
with
\[\alpha(\rho_0)=\int_0^1\frac{\partial H^{N}_{0}(\rho_0, t)}{\partial \rho}\,d t,\quad |\partial_{\rho_0}\, \alpha(\rho_0)|\ge C\, A^n>0,\]
and
$$|\partial^{p}_{\rho_{0}}\partial^{q}_{\xi_{0}}F|\leq C\varepsilon_{0},\;\;
|\partial^{p}_{\rho_{0}}\partial^{q}_{\xi_{0}}G|\leq C\varepsilon_{0},\;\;p+q\leq 5.$$
Since \eqref{eq72} is Hamiltonian, the map $P$ is symplectic.
By Moser's twist theorem at {pp.50-54} of \cite{Moser}, $\mathcal{P}$ has an invariant curve $\Gamma$ in the annulus  $[2,3]\times \mathbb{T}^1.$
Since $A$ can be arbitrarily large, it follows that the time-1 map of the original system has  an invariant curve $\Gamma_A$ in the annulus  $[2\,A+C ,3\,A-C]\times \mathbb{T}^1$ with $C$ is a constant independent of $A$. Choosing a sequence $A=A_k\to\infty$ as $k\to\infty$, we have that there are
countable many invariant curves $\Gamma_{A_k}$, clustering at $\infty.$ Therefore any solution of the original system is bounded. This completes the proof of Theorem.

\begin{rem} Any solutions starting from the invariant curves $\Gamma_{A_k}$ ($k=1,2,...$) are quasi-periodic with frequencies $(1,\omega_k)$ in time $t$, where $(1,\omega_k)$ satisfies Diophantine conditions and $\omega>C\, A_k^n$.

\end{rem}

\end{document}